\documentclass[11pt]{amsart}

\usepackage{amsmath,amsthm,amsfonts,amscd,flafter,epsf,amssymb}
\usepackage{epsfig, psfrag, subfigure, bm}
\usepackage[margin=1in]{geometry}

\newtheorem{thm}{Theorem}[section]
\newtheorem{defn}[thm]{Definition}
\newtheorem{lemma}[thm]{Lemma}
\newtheorem{prop}[thm]{Proposition}

\newcommand{\calh}{\mathcal{H}}

\begin{document}

\title[{An algorithm for computing some Heegaard Floer homologies}]{An algorithm for computing some Heegaard Floer homologies}

\author[Sucharit Sarkar]{Sucharit Sarkar}
\address{Department of Mathematics, Princeton University, Princeton, NJ 08544}
\email{sucharit@math.princeton.edu}

\author[Jiajun Wang]{Jiajun Wang}
\address{Department of Mathematics, California Institute of Technology, Pasadena, CA 91125}
\email{wjiajun@caltech.edu}

\begin{abstract}
In this paper, we give an algorithm to compute the hat version of Heegaard Floer homology of a closed oriented three-manifold. This method also allows us to compute the filtration coming from a null-homologous link in a three-manifold.
\end{abstract}


\date{May 1st, 2008}

\maketitle

\section{Introduction}

Heegaard Floer homology is a collection of invariants for closed oriented three-manifolds, introduced by Peter Ozsv\'{a}th and Zolt\'{a}n Szab\'{o} \cite{OSzClosed3M, OSzClosed3MApp}. There are four versions, denoted by $\widehat{HF}, HF^{\infty}, HF^{+}$ and $HF^{-}$, which are graded abelian groups. The hat version $\widehat{HF}(Y)$ is defined as the homology of a chain complex $\widehat{CF}(Y)$ coming from a Heegaard diagram of the three-manifold $Y$. The differentials count the number of points in certain moduli spaces of holomorphic disks, which are hard to compute in general.

There is also a relative version of the theory corresponding to pairs $(Y,K)$, where $K$ is a knot in $Y$. If $K$ is null-homologous, then a Seifert surface $S$ of $K$ induces a filtration of the chain complex $\widehat{CF}(Y)$, and the chain homotopy type of the filtered chain complex is a knot invariant. The homology groups $\widehat{HFK}(Y,K)$ of successive quotients of filtration levels are called knot Floer homology groups (\cite{OSzKnotInv, Rasmussen, OSzLinkInv}).

A cobordism between two three-manifolds induces homomorphisms on the Heegaard Floer homology groups of the two three-manifolds. In fact, the homomorphisms on $HF^{-}$ and $HF^{+}$ can be used to construct an invariant of smooth four-manifolds with $b_2^+>1$ (\cite{OSzSmooth4M}), called the Ozsv\'ath-Szab\'o invariant. Conjecturally, the Ozsv\'ath-Szab\'o four-manifold invariant is equivalent to the gauge-theoretic Seiberg-Witten invariant.

Heegaard Floer homology turns out to be a fruitful and powerful theory in the study of three-dimensional and four-dimensional topology. It gives an alternate proof of the Donaldson diagonalization theorem and the Thom conjecture for $\mathbb{CP}^2$ (\cite{OSzAbsoluteGrade}). Heegaard Floer homology also detects the Thurston norm of a three-manifold (\cite{OSzThurstonNorm, N522}). Moreover, knot Floer homology detects the genus (\cite{OSzGenusBound}) and fiberedness (\cite{Ghiggini, NiFibred, Juhasz}) of knots and links in the three-sphere. There is an invariant $\tau$  coming from the knot filtration, whose absolute value gives a lower-bound of the slice genus for knots in the three-sphere (\cite{OSz4BallGenus}). 

Despite its success, there was no general method to compute the invariants. There were combinatorial descriptions in certain special cases, but the computation for an arbitrary three-manifold was an open problem (\cite{Differentfaces}). In this paper, we give an algorithm to compute $\widehat{HF}(Y)$ for a three-manifold $Y$, and also $\widehat{HFK}(Y,K)$ for a knot $K$ in any three-manifold. All our computations will be done with coefficients in $\mathbb{F}_2=\mathbb{Z}/2\mathbb{Z}$. We show that one can always find Heegaard diagrams satisfying certain properties (Definition \ref{defn:nicediagram}). Using such Heegaard diagrams, which we call nice, it will be easy to compute $\widehat{HF}$ and $\widehat{HFK}$. Our main results are summarized in the following theorems.

\begin{thm}\label{thm:holodisk} Given a nice Heegaard diagram of a closed oriented three-manifold $Y$, $\widehat{HF}(Y)$ can be computed combinatorially. Similarly, for a knot $K\subset Y$, $\widehat{HFK}(Y,K)$ can be computed combinatorially in a nice Heegaard diagram.
\end{thm}

\begin{thm}\label{thm:nicediagram} Every closed oriented three-manifold $Y$ admits a nice Heegaard diagram. For a null-homologous knot $K$ in a closed oriented three-manifold $Y$, the pair $(Y,K)$ admits a compatible nice Heegaard diagram. In fact, there is an algorithm to convert any pointed Heegaard diagram to a nice Heegaard diagram via isotopies and handleslides.
\end{thm}

It will be interesting to compare our result with the recent work of Ciprian Manolescu, Peter Ozsv{\'a}th and the first author in \cite{ManOzsSar}, where they gave a combinatorial description of knot Floer homology of knots in $S^3$, in all versions. 

We hope this method can be generalized to compute some of the other versions, notably $HF^{-}(Y)$ and $HF^{+}(Y)$. It would also be nice to have a  proof of the invariance of the combinatorial description without using holomorphic disks.

The paper is organized as follows. In Section \ref{section-prelims}, we give an overview of certain concepts in Heegaard Floer theory. In Section \ref{section-holo}, we give a combinatorial characterization of index one holomorphic disks in nice Heegaard diagrams. In Section \ref{section-algorithm}, we give an algorithm to get such Heegaard diagrams. In Section \ref{section-examples}, we give examples to demonstrate our algorithm for three-manifolds and knots in the three-sphere.

\section*{Acknowledgement}
The first author wishes to thank Zolt\'{a}n Szab\'{o} for introducing him to the subject of Heegaard Floer homology and having many helpful discussions at various points. 

This work was done when the second author was an exchange graduate student in Columbia University. He is grateful to the Columbia math department for its hospitality. He also would like to thank Rob Kirby and Peter Ozsv{\'a}th for their continuous guidance, support and encouragement.

We thank Matthew Hedden, Robert Lipshitz, Ciprian Manolescu, Peter Ozsv{\'a}th, Jacob Rasmussen and Dylan Thurston for making comments and having helpful discussions during the development of this work.

We also thank the referees for their comments and suggestions.

\section{Preliminaries}\label{section-prelims}

In this section, we review the definition of Heegaard Floer homology. See \cite{OSzClosed3M, OSzClosed3MApp, OSzLinkInv, Lipshitz} for details.

\subsection{Definition of $\widehat{HF}$} 

The Heegaard Floer homology of a closed oriented three-manifold $Y$ is defined from a pointed Heegaard diagram representing $Y$. 

A Heegaard splitting of $Y$ is a decomposition of $Y$ into two handlebodies glued along their boundaries. We fix a self-indexing Morse function $f$ on $Y$ with $k$ index zero critical points and $k$ index three critical points. (We usually choose $k=1$.) Then $f$ gives a Heegaard splitting of $Y$, where the two handlebodies are given by $f^{-1}(-\infty, \frac{3}{2}]$ and $f^{-1}[\frac{3}{2},\infty)$. If the number of index one critical points or the number of index two critical points of $f$ is $(g+k-1)$, then $\Sigma=f^{-1}({3}/{2})$ is a genus $g$ surface. We fix a gradient like flow on $Y$ corresponding to $f$. We require $f$ to have the property that $Y$ contains a disjoint union of $k$ flow lines, each flowing from an index zero critical point to an index three critical point. We get a collection $\bm{\alpha}=(\alpha_1,\cdots,\alpha_{g+k-1})$ of $\alpha$ circles on $\Sigma$ which flow down to the index one critical points, and another collection $\bm{\beta}=(\beta_1,\cdots,\beta_{g+k-1})$ of $\beta$ circles  on $\Sigma$ which flow up to the index two critical points. Note that both $\Sigma\setminus \bm{\alpha}$ and $\Sigma\setminus \bm{\beta}$ have $k$ components. 

We fix $k$ points $w_1,\ldots,w_k$ (called {\it basepoints}) in the complement of the $\alpha$ circles and the $\beta$ circles in $\Sigma$, such that each component of $\Sigma \setminus \bm{\alpha}$ contains exactly one $w_i$ and each component of $\Sigma\setminus \bm{\beta}$ contains exactly one $w_j$. This is equivalent to the condition that the trajectories of $w_i$'s under the gradient like flow, hit all the index zero and all the index three critical points. We write $\bm w=(w_1,\cdots,w_k)$. The tuple $(\Sigma, \bm{\alpha},\bm{\beta}, \bm{w})$ is called a {\it pointed Heegaard diagram} for $Y$.

There are some moves on a Heegaard diagram that do not change the underlying three-manifold. An {\it isotopy} moves the $\alpha$ curves and $\beta$ curves in two one-parameter families $\bm\alpha_t$ and $\bm\beta_t$ in $\Sigma\setminus\bm w$, moving by isotopy, such that the $\alpha$ curves remain disjoint and the $\beta$ curves remain disjoint for each $t$. In a {\it handleslide} of $\bm{\alpha}$, we replace a pair of $\alpha$ curves $\alpha_i$ and $\alpha_j$ with a pair $\alpha_i$ and $\alpha^{\prime}_j$, such that the three curves $\alpha_i$, $\alpha_j$ and $\alpha^{\prime}_j$ bound a pair of pants in $\Sigma\setminus\bm w$ disjoint from all the other $\alpha$ curves. A handleslide of $\bm{\beta}$ is defined similarly. There is also a move called {\it stabilization}, but we will not be using it in the present paper. These moves are called {\it Heegaard moves}.


Heegaard Floer homology is a certain version of Lagrangian Floer homology. The ambient symplectic manifold is the symmetric product ${\rm Sym}^{g+k-1}(\Sigma)$. The two half-dimensional totally real subspaces are the tori $T_{\alpha}=\alpha_1\times\cdots\times\alpha_{g+k-1}$ and $T_{\beta}=\beta_1\times\cdots\times\beta_{g+k-1}$. The generators for the chain complex $\widehat{CF}(\Sigma,\bm\alpha,\bm\beta,\bm w)$ are the intersection points between these two tori, and the boundary maps are given by counting certain holomorphic disks. For more details see \cite{OSzClosed3M}, and see \cite{OSzLinkInv} for the issue of {\it boundary degenerations} when $k>1$.

\begin{thm}[Ozsv\'ath-Szab\'o, \cite{OSzClosed3M, OSzLinkInv}] When $k=1$, the homology of the chain complex $\widehat{CF}(\Sigma,\bm\alpha,\bm\beta,\bm w)$ is an invariant for the three-manifold $Y$, written as $\widehat{HF}(Y)$. For a general $k$, and $Y$ a rational homology three-sphere, we have 
$${H_*(\widehat{CF}(\Sigma,\bm\alpha,\bm\beta,\bm w))}\cong \widehat{HF}(Y)\otimes H_*(T^{k-1}),$$
where $H_*(T^{k-1})$ is the singular homology of the
$(k-1)$-dimensional torus with coefficients in $\mathbb{F}_2$.
\end{thm}

When we have a link $L$ in $Y$, we ensure that $L$ is a union of flow lines from index zero critical points to index three critical points. We also ensure that $L$ contains all index zero and index three critical points and contains no index one or two critical points. We orient $L$ and $\Sigma$ and define $w_i$'s as the positive intersection points between $L$ and $\Sigma$. We write the other $k$ intersections as $\bm z=(z_1,\cdots, z_k)$. Such a Heegaard diagram,  denoted by $(\Sigma,\bm\alpha,\bm\beta,\bm w,\bm z)$, is called a {\it pointed Heegaard diagram} for the pair $(Y,L)$. (For a link with $l$ components, we usually choose $k=l$.) The knot (link) Floer homology $\widehat{HFK}(Y,L)$ is defined similarly, where the boundary maps count a more restricted class of holomorphic disks. See \cite{OSzKnotInv, Rasmussen, OSzLinkInv} for details.

\subsection{Cylindrical reformulation of $\widehat{HF}$} In the present paper, we will use the cylindrical reformulation of the Heegaard Floer homology by Lipshitz. See \cite{Lipshitz} for details.

Given a pointed Heegaard diagram $(\Sigma,\bm\alpha,\bm\beta,\bm w)$, the generators of the chain complex $\widehat{CF}$ are given by formal sums of $(g+k-1)$ distinct points in $\Sigma$, $\bm x=x_1+\cdots+x_{g+k-1}$, such that each $\alpha$ circle contains some $x_i$ and each $\beta$ circle contains some $x_j$. A connected component of $\Sigma\setminus(\bm\alpha\cup\bm\beta)$ is called a {\it region}. A formal sum of regions with integer coefficients is called a {\it $2$-chain}. Given two generators $\bm x$ and $\bm y$, we define $\pi_2(\bm x,\bm y)$ to be the collection of all $2$-chains $\phi$ such that $\partial(\partial(\phi)|_{\alpha})=\bm y-\bm x$. Such $2$-chains are called {\it domains}. Given a point $p\in \Sigma\setminus(\bm\alpha\cup\bm\beta)$, let $n_p(\phi)$ be the coefficient of the region containing $p$ in $\phi$. A domain $\phi$ is {\it positive} if $n_p(\phi)\geq 0$ for all points $p\in\Sigma\setminus(\bm\alpha\cup\bm\beta)$.  We define $\pi_2^0(\bm x,\bm y)=\{\phi\in\pi_2(\bm x,\bm y)\ |\ n_{w_i}(\phi)=0\ \forall i\}$. A Heegaard diagram is {\it admissible}, if, for every generator $\bm x$, any positive domain $\phi\in \pi_2^0(\bm x,\bm x)$ is trivial. If the three-manifold $Y$ has $b_1(Y)>0$, we require the Heegaard diagram to be admissible.

Fix two generators $\bm x$, $\bm y$ and a domain $\phi\in\pi_2^0(\bm x,\bm y)$. Let $S$ be a surface with boundary, with $2(g+k-1)$ marked points $(X_1,\cdots,X_{g+k-1},Y_1,\cdots,Y_{g+k-1})$ on $\partial S$, such that the $X$ points and the $Y$ points alternate. The $2(g+k-1)$ arcs on $\partial S$ in the complement of the marked points are divided into two groups $A$ and $B$, each containing $(g+k-1)$ arcs, such that the $A$ arcs and the $B$ arcs alternate. Let $p_1$ and $p_2$ be the projection maps from $\Sigma\times D^2$ onto its first and second factors. Look at maps $u:S\rightarrow\Sigma\times D^2$ such that the image of $p_1\circ u$ is $\phi$ (as $2$-chains) and the image of $p_2\circ u$ is $(g+k-1)D^2$ (in second homology). We also want the $X$ points on $\partial S$ to map injectively by $p_1\circ u$ to the $x_i$'s and to map to $-i$ in the unit disk by $p_2\circ u$. Similarly we want the $Y$ points to map injectively to the $y_i$'s by $p_1\circ u$ and to $i$ in the unit disk by $p_2\circ u$. Furthermore we also require the $A$ arcs in $\partial S$ to map to $\alpha$ arcs by $p_1\circ u$ and under $p_2\circ u$ to map to the arc $e_1$ in $\partial(D^2)$ joining $-i$ to $i$ in half-plane $Re(s)>0$. Similarly we require the $B$ arcs to map to $\beta$ arcs by $p_1\circ u$ and to map to the arc $e_2$ in $\partial(D^2)$ in the half-plane $Re(s)<0$ by $p_2\circ u$.

Now fix complex structures on $\Sigma$ and $D^2$ and take the product complex structure on $\Sigma\times D^2$. A generic perturbation gives an almost complex structure which achieves transversality for the homology class $\phi$. In our case, we can achieve this by a generic perturbation of the $\alpha$ curves and the $\beta$ curves (\cite[Lemma 3.10]{Lipshitz}).

The holomorphic embeddings $u$ which satisfy the above conditions and whose homology class is $\phi$ form a moduli space, which we denote by ${\mathcal M}(\phi)$. The {\it Maslov index} $\mu(\phi)$ of $\phi$ gives the expected dimension of ${\mathcal M}(\phi)$. It can be computed combinatorially in terms of the Euler measure and the point measures, which are defined as follows. For a generator $\bm x=\sum x_i$ and a domain $\phi$, $\mu_{x_i}(\phi)$ is defined to be the average of the coefficients of the four regions around $x_i$ in $\phi$. The {\it point measure} $\mu_{\bm x}(\phi)$ is defined as $\sum \mu_{x_i}(\phi)$. If we fix a metric on $\Sigma$ which makes all the $\alpha$ and $\beta$ circles geodesic, intersecting each other with right angles, then the {\it Euler measure} $e(\phi)$ is defined to be $\frac{1}{2\pi}$ of the integral of the curvature on $\phi$. The Euler measure is clearly additive, and if $D$ is a $2n$-gon region, then $e(D)=1-\frac{n}{2}$.

\begin{prop}[Lipshitz, \cite{Lipshitz}] For a domain $\phi\in\pi_2(\bm x,\bm y)$, the Maslov index is given by
\begin{eqnarray*}
\mu(\phi)=e(\phi)+\mu_{\bm x}(\phi)+\mu_{\bm y}(\phi)
\end{eqnarray*}
\end{prop}

If $\phi$ is non-trivial, the moduli space ${\mathcal M}(\phi)$ admits a free $\mathbb R$-action coming from the one-parameter family of holomorphic automorphisms of $D^2$ which preserve $\pm i$ and the boundary arcs $e_1$ and $e_2$. In particular, if $\mu(\phi)=1$, the {\it unparametrized moduli space} ${\mathcal M}(\phi)/{\mathbb R}$ is a zero-dimensional manifold, and then the {\it count function} $c(\phi)$ is defined to be the number of points in ${\mathcal M}(\phi)/{\mathbb R}$, counted modulo $2$. The boundary map in the chain complex $\widehat{CF}$ is given by
\begin{eqnarray*}
\partial \bm x=\sum_{\bm y}\sum_{\{\phi\in\pi_2^0(\bm x,\bm y)\ |\ \mu(\phi)=1\}}c(\phi)\bm y.
\end{eqnarray*}

\begin{thm}[Lipshitz, \cite{Lipshitz}] For a three-manifold $Y$, the homology of the chain complex $(\widehat{CF},\partial)$ is isomorphic to $H_*(\widehat{CF}(\Sigma, \bm\alpha,\bm\beta,\bm w))$.\end{thm} 

Note that the only non-combinatorial part of the theory is the count function $c(\phi)$. 

\subsection{Positivity of domains with holomorphic representatives} We will need the following proposition, which asserts that only positive domains can have holomorphic representatives.

\begin{prop} Let $\phi$ be a domain in $\pi_2^0(\bm x,\bm y)$. If $\phi$ has a holomorphic representative, then $\phi$ is a positive domain. In particular, if $c(\phi)\neq 0$, then $\phi$ is a positive domain.
\end{prop}

\begin{proof}
If $\phi$ has a holomorphic representative, then there exists some holomorphic embedding $u$ of the type described above. Then for any point $p\in\Sigma\setminus(\bm\alpha\cup\bm\beta)$, $n_p(\phi)$ is simply the intersection number of $u(S)$ and $\{p\}\times D^2$. Since both of them are holomorphic objects in the product complex structure, they have positive intersection number and hence $n_p(\phi)\ge 0$. Here we require the complex structure on $\Sigma\times D^2$ to be standard near the basepoints. See \cite{Lipshitz} for a general discussion.
\end{proof}

If a domain $\phi$ has a holomorphic representative, the number of branch points of $p_2\circ u$ is given by $\mu_{\bm x}+\mu_{\bm y}-e(\phi)$ (\cite{Lipshitz, Rasmussen}). Furthermore, in such a situation the Maslov index can also be calculated as $\mu(\phi)=2e(\phi)+g+k-1-\chi(S)=e(\phi)+b+\frac{1}{2}(g+k-1-t)$, where $b$ denotes the number of branch points of $p_1\circ u$, and $t$ denotes the number of {\it trivial disks}, i.e. the components of $S$ which are mapped to a point by $p_1\circ u$ (which correspond to coordinates $x_i$ of $\bm x$ with $\mu_{x_i}=0$).  

\section{Holomorphic disks in nice Heegaard diagrams}\label{section-holo}

In this section, we study index one holomorphic disks in nice Heegaard diagrams.

\begin{defn}\label{defn:nicediagram} Let $\calh=(\Sigma, \bm\alpha, \bm\beta, \bm w)$ be a pointed Heegaard diagram for a three-manifold $Y$. $\calh$ is called {\bf nice} if any region that does not contain any basepoint $w_i$ in $\bm w$ is either a bigon or square.\end{defn}

Let $Y$ be a closed oriented three-manifold. Suppose $Y$ has a nice admissible Heegaard diagram $\calh=(\Sigma, \bm\alpha, \bm\beta, \bm w)$. We choose a product complex structure on $\Sigma\times D^2$. 

\begin{defn}\label{defn:empty}
A domain $\phi\in \pi_2^0(\bm x,\bm y)$ with coefficients $0$ and $1$ is called an {\bf empty embedded $2n$-gon}, if it is topologically an embedded disk with $2n$ vertices on its boundary, such that at each vertex $v$, $\mu_v(\phi)=\frac{1}{4}$, and it does not contain any $x_i$ or $y_i$ in its interior.
\end{defn}

The following two theorems show that, for a domain $\phi\in\pi_2^0(\bm x,\bm y)$, the count function $c(\phi)\neq 0$ if and only if $\phi$ is an empty embedded bigon or an empty embedded square, and in that case $c(\phi)=1$. Thus $c(\phi)$ can be computed combinatorially in a nice Heegaard diagram.

\begin{thm}\label{thm:holo=>embedded}
Let $\phi\in \pi_2^0(\bm x,\bm y)$ be a domain such that $\mu(\phi)=1$. If $\phi$ has a holomorphic representative, then $\phi$ is an empty embedded bigon or an empty embedded square.
\end{thm}

\begin{proof}

We know that only positive domains can have holomorphic representatives. We also know that bigons and squares have non-negative Euler measure. We will use these facts to limit the number of possible cases.

Suppose $\phi=\sum a_i D_i$, where $D_i$'s are regions containing no
basepoints. Since $\phi$ has a holomorphic representative, we have
$a_i\ge 0, \forall i$. Since each $D_i$ is a bigon or a square, we have
$e(D_i)\ge 0$ and hence $e(\phi)\geq0$. So, by Lipshitz' formula
$\mu(\phi)=e(\phi)+\mu_{\bm x}(\phi)+\mu_{\bm y}(\phi)$, we get $0\le
\mu_{\bm x} +\mu_{\bm y}\le 1$ .

Now let $\bm x=x_1+\cdots+x_g$ and $\bm y=y_1+\cdots+y_g$, with $x_i,y_i \in \alpha_i$. We say $\phi$ hits some $\alpha$ circle if $\partial \phi$ is non-zero on some part of that $\alpha$ circle. Since $\phi \neq n\Sigma$, it has to hit at least one $\alpha$ circle, say $\alpha_1$, and hence $\mu_{x_1}, \mu_{y_1} \ge \frac{1}{4}$ as $\partial(\partial \phi_{|\alpha})=\bm y-\bm x$. Also if $\phi$ does not hit $\alpha_i$, then $x_i=y_i$ and they must lie outside the domain $\phi$, since otherwise we have $\mu_{x_i}=\mu_{y_i}\geq \frac{1}{2}$ and hence $\mu_{\bm x}+\mu_{\bm y}$ becomes too large.

We now note that $e(\phi)$ can only take half-integral values, and thus only the following cases might occur.

\begin{itemize}

\item {\bf Case 1.} $\phi$ hits $\alpha_1$ and another $\alpha$ circle, say $\alpha_2$, $\phi$ consists of squares, $\mu_{x_1}=\mu_{x_2}=\mu_{y_1}=\mu_{y_2}=\frac{1}{4}$, and there are $(g+k-3)$ trivial disks.

\item {\bf Case 2.} $\phi$ hits $\alpha_1$, $D(\phi)$ consists of squares and exactly one bigon, $\mu_{x_1}=\mu_{y_1}=\frac{1}{4}$, and there are $(g+k-2)$ trivial disks.

\item {\bf Case 3.} $\phi$ hits $\alpha_1$, $D(\phi)$ consists of squares, $\mu_{x_1}+\mu_{y_1}=1$, and there are $(g+k-2)$ trivial disks.

\end{itemize}

Using the reformulation by Lipshitz, in each of these cases, we will try to figure out the surface $S$ which maps to $\Sigma\times D^2$. Recall that a trivial disk is a component of $S$ which maps to a point in $\Sigma$ after post-composing with the projection $\Sigma\times D^2\rightarrow\Sigma$.

The first case corresponds to a map from $S$ to $\Sigma$ with $\chi(S)=(g+k-2)$, and $S$ has $(g+k-3)$ trivial disk components. If the rest of $S$ is $F$, then $F$ is a double branched cover over $D^2$ with $\chi(F)=1$ and $1$ branch point (for holomorphic maps, the number of branch points is given by $\mu_{\bm x} + \mu_{\bm y} - e(\phi)$), i.e. $F$ is a disk with $4$ marked points on its boundary. Call the marked points corners, and call $F$ a square.

In the other two cases, $S$ has $(g+k-2)$ trivial disk components, so if $F$ denotes the rest of $S$, then $F$ is just a single cover over $D^2$. Thus the number of branch points has to be $0$. But in the third case the number of branch points is $1$, so the third case cannot occur. In the second case, $F$ is a disk with $2$ marked points on its boundary. Call the marked points corners, and call $F$ a bigon. 

Thus in both the first and the second cases, $\phi$ is the image of $F$ and all the trivial disks map to the $x$-coordinates (which are also the $y$-coordinates) which do not lie in $\phi$. Note that in both cases, the map from $F$ to $\phi$ has no branch point, so it is a local diffeomorphism, even at the boundary of $F$. Furthermore using the condition that $\mu_{x_i}$ (or $\mu_{y_i}$) $=\frac{1}{4}$ whenever it is non-zero, we conclude that there is exactly one preimage for the image of each corner of $F$.

All we need to show is that the map from $F$ to $\Sigma$ is an embedding, or in other words, the local diffeomorphism from $F$ to $\phi$ is actually a diffeomorphism. We will prove this case by case.

\subsection*{Case 1}

In this case we have an immersion $f:F\rightarrow \Sigma$, where $F$ is a square (with boundary). Look at the preimage of all the $\alpha$ and $\beta$ circles in $F$. Using the fact $f$ is a local diffeomorphism, we see that each of the preimages of $\alpha$ and $\beta$ arcs are also $1$-manifolds, and by an abuse of notation, we will also call them $\alpha$ or $\beta$ arcs. Using the embedding condition near the $4$ corners, we see that at each corner only one $\alpha$ arc and only one $\beta$ arc can come in. The different $\alpha$ arcs cannot intersect and the different $\beta$ arcs cannot intersect, and all intersections between $\alpha$ and $\beta$ arcs are transverse.

Note that since the preimage of each square region is a square, $F$ (with all the $\alpha$ and $\beta$ arcs) is also tiled by squares. Thus the $\alpha$ arcs in $F$ cannot form a closed loop, for in that case $F\setminus$\{inside of loop\} has negative Euler measure and hence cannot be tiled by squares. Similarly the $\beta$ arcs cannot form a loop. Also no $\alpha$ arc can enter and leave $F$ through the same $\beta$ arc on the boundary, for again the outside will have negative Euler measure. Thus the $\alpha$ arcs slice up $F$ into vertical rectangles, and in each rectangle, no $\beta$ arc can enter and leave through the same $\alpha$ arc. This shows that the $\alpha$ arcs and $\beta$ arcs make the standard co-ordinate chart on $F$, as in Figure \ref{fig:holo_square}.
\begin{figure}[htbp]
 \center{\includegraphics{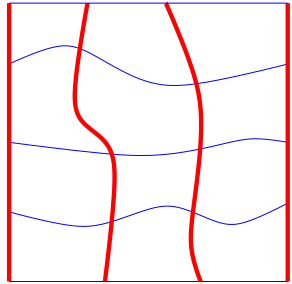}}
 \caption{\label{fig:holo_square}
  {\bf Preimage of $\alpha$ and $\beta$ arcs for a square.} We make the convention for all figures in the paper that the thick solid arcs denote $\alpha$ arcs and the thin solid arcs denote $\beta$ arcs.
 }
\end{figure}

We call the intersection points between $\alpha$ and $\beta$ arcs in $F$ vertices (and we are still calling the four original vertices on the boundary of the square $F$ corners). Note that to show $f$ is an embedding, it is enough to show that no two different vertices map to the same point. Assume $p,q\in F$ are distinct vertices with $f(p)=f(q)$. There could be two subcases.
\begin{itemize}
 \item Both $p$ and $q$ are in $\mathring{F}$.
 \item At least one of $p$ and $q$ is in $\partial F$.
\end{itemize}

We will reduce the first subcase to the second. Assume both $p$ and $q$ are in the interior of $F$. Choose a direction on the $\alpha$ arc passing through $f(p)=f(q)$ in $\Sigma$, and keep looking at successive points of intersection with $\beta$ arcs, and locate their inverse images in $F$. For each point, we will get at least a pair of inverse images, one on the $\alpha$ arc through $p$, and one on the $\alpha$ arc through $q$, until one of the points falls on $\partial F$, and thus we have reduced it to the second subcase.

In the second subcase, without loss of generality, we assume that $p$ lies on a $\beta$ arc on $\partial F$. Then choose a direction on the $\beta$ arc in $\Sigma$ through $f(p)=f(q)$ and proceed as above, until one of the preimages hits an $\alpha$ arc on $\partial F$. If that preimage is on the $\beta$ arc through $q$, then reverse the direction and proceed again, and this time we can ensure that the preimage which hits $\alpha$ arc on $\partial F$ first is the one that was on the $\beta$ arc through $p$. Thus we get $2$ distinct vertices in $F$ mapping to the same point in $\Sigma$, one of them being a corner. This is a contradiction to the embedding assumption near the corners.

\subsection*{Case 2}

In this case we have an immersion $f:F\rightarrow \Sigma$ with $F$ being a bigon. Again look at the preimage of $\alpha$ and $\beta$ circles. All intersections will be transverse (call them vertices), and at each of the $2$ corners there can be only one $\alpha$ arc and only one $\beta$ arc. Again there cannot be any closed loops. We get an induced tiling on $F$ with squares and $1$ bigon.

This time the $\alpha$ arcs can (in fact they have to) enter and leave $F$ through the same $\beta$ arc, but they have to do it in a completely nested fashion, i.e., there is only one bigon piece in $F\setminus\bm\alpha$, the ``innermost bigon". Thus $F$ decomposes into two pieces, the innermost bigon and the rest. In case there are no $\alpha$ arcs in $\mathring{F}$, the rest might be empty, but otherwise it is a square. From the arguments in the earlier case, the $\beta$ arcs must cut up the square piece in a standard way, and from the previous argument the $\beta$ arcs must enter and leave the bigon in a nested fashion, as in Figure \ref{fig:holo_bigon}.
\begin{figure}[htbp]
 \center{\includegraphics{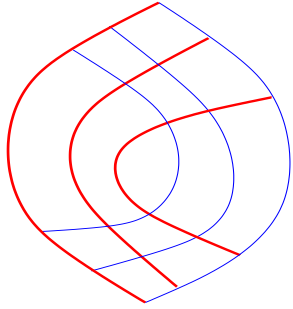}}
 \caption{\label{fig:holo_bigon}
  {\bf Preimage of $\alpha$ and $\beta$ arcs for a bigon.} 
 }
\end{figure}

Again to show $f$ is an embedding, it is enough to show that it is an embedding restricted to vertices. Take $2$ distinct vertices $p,q$ mapping to the same point, and follow them along $\alpha$ arcs in some direction, until one of them hits a $\beta$ arc on $\partial F$. Then follow them along $\beta$ arcs, and there exists some direction such that one of them will actually hit a corner, giving the required contradiction.

So in either case, $f$ is an embedding.
\end{proof}

\begin{thm}\label{thm:embedded=>holo}
If $\phi\in \pi_2^0(\bm x,\bm y)$ is an empty embedded bigon or an empty embedded square, then the product complex structure on $\Sigma\times D^2$ achieves transversality for $\phi$ under a generic perturbation of the $\alpha$ and the $\beta$ circles, and $\mu(\phi)=c(\phi)=1$.
\end{thm}

\begin{proof}

Let $\phi$ be an empty embedded $2n$-gon. Each of the corners of $\phi$ must be an $x$-coordinate or a $y$-coordinate, and at every other $x$ (resp. $y$) coordinate the point measure $\mu_{x_i}$ (resp. $\mu_{y_i}$) is zero. Therefore $\mu_{\bm x}(\phi)+\mu_{\bm y}(\phi)=2n\cdot\frac{1}{4}=\frac{n}{2}$. Also $\phi$ is topologically a disk, so it has Euler characteristic $1$. Since it has $2n$ corners each with an angle of $\frac{\pi}{4}$, the Euler measure $e(\phi)=1-\frac{2n}{4}=1-\frac{n}{2}$. Thus the Maslov index $\mu(\phi)=1$.

By \cite[Lemma 3.10]{Lipshitz}, we see that $\phi$ satisfies the boundary injective condition, and hence under a generic perturbation of the $\alpha$ and the $\beta$ circles, the product complex structure achieves transversality for $\phi$.

When $\phi$ is an empty embedded square, we can choose $F$ to be a disk with $4$ marked points on its boundary, which is mapped to $\phi$ diffeomorphically. Given a complex structure on $\Sigma$, the holomorphic structure on $F$ is determined by the cross-ratio of the four points on its boundary, and there is an one-parameter family of positions of the branch point in $D^2$ which gives that cross-ratio. Thus there is a holomorphic branched cover $F\rightarrow D^2$ satisfying the boundary conditions, unique up to reparametrization. Hence $\phi$ has a holomorphic representative, and from the proof of Theorem \ref{thm:holo=>embedded} we see that this determines the topological type of $F$, and hence it is the unique holomorphic representative.

When $\phi$ is an empty embedded bigon, we can choose $F$ to be a disk with $2$ marked points on its boundary, which is mapped to $\phi$ diffeomorphically. A complex structure on $\Sigma$ induces a complex structure on $F$, and there is a unique holomorphic map from $F$ to the standard $D^2$ after reparametrization. Thus again $\phi$ has a holomorphic representative, and similarly it must be the unique one.
\end{proof}

\begin{proof}[\bf Proof of Theorem \ref{thm:holodisk}]

Theorems \ref{thm:holo=>embedded} and \ref{thm:embedded=>holo} make the count function $c(\phi)$ combinatorial in a nice Heegaard diagram. For a domain $\phi\in \pi_2^0(\bm x,\bm y)$ with $\mu(\phi)=1$, we have $c(\phi)=1$ if $\phi$ is an empty embedded bigon or an empty embedded square, and $c(\phi)=0$ otherwise.
\end{proof}

\newpage

\section{Algorithm to get nice Heegaard diagrams}\label{section-algorithm}

In this section, we prove Theorem \ref{thm:nicediagram}. We will demonstrate an algorithm which, starting with an admissible pointed Heegaard diagram, gives an admissible nice Heegaard diagram by doing isotopies and handleslides on the $\beta$ curves.

For a Heegaard diagram, we call bigon and square regions {\it good} and all other regions {\it bad}. We will first do some isotopies to ensure all the regions are disks. We will then define a complexity for the Heegaard diagram which attains its minimum only if all the regions not containing the basepoints are good. We will do an isotopy or a handleslide which will decrease the complexity if the complexity is not the minimal one.

\subsection{The algorithm}

Let $\calh=(\Sigma, \bm\alpha,\bm\beta, w)$ be a pointed Heegaard diagram with a single basepoint $w$. We consider Heegaard diagrams with more basepoints in the last subsection.

\subsection*{Step 1. Killing non-disk regions.}\

We do finger moves on $\beta$ circles to create new intersections with $\alpha$ circles. After doing this sufficiently many times, every region in $\calh$ becomes a disk. We first ensure that every $\alpha$ circle intersects some $\beta$ circle and every $\beta$ circle intersects some $\alpha$ circle.

\begin{figure}[htbp]
 \center{\includegraphics[width=350pt]{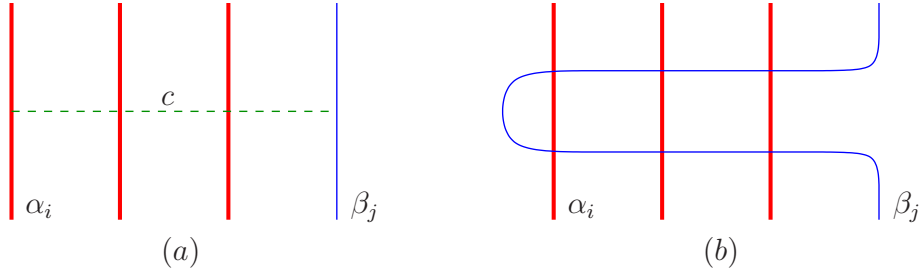}}
 \vspace{-10pt}
 \caption{\label{fig:no_isolate_alpha}
  {\bf Making each $\alpha$ circle intersect some $\beta$ circle.}
 }
\end{figure}

If $\alpha_i$ does not intersect any $\beta$ circle, we can find an arc $c$ connecting $\alpha_i$ to some $\beta_j$ avoiding the intersections of $\alpha$ and $\beta$ circles, as indicated in Figure \ref{fig:no_isolate_alpha}(a). We can select $c$ such that $c$ intersects $\bm\beta$ just at the endpoint. Doing a finger move of $\beta_j$ along $c$ as in Figure \ref{fig:no_isolate_alpha}(b) will make $\alpha_i$ intersect some $\beta$ circle.

\begin{figure}[htbp]
 \center{\includegraphics[width=350pt]{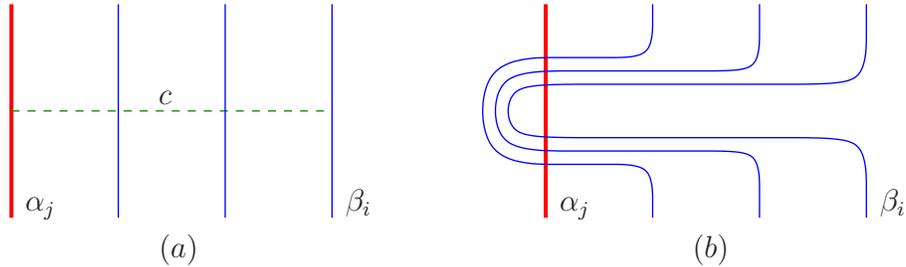}}
 \vspace{-10pt}
 \caption{\label{fig:no_isolate_beta}
  {\bf Making each $\beta$ circle intersect some $\alpha$ circle.}
 }
\end{figure}

Similarly, if $\beta_i$ does not intersect any $\alpha$ circle, we find an arc $c$ connecting $\beta_i$ to some $\alpha_j$ so that $c\cap\bm\alpha$ contains a single point as in Figure \ref{fig:no_isolate_beta}(a). We then do the operation as depicted in Figure \ref{fig:no_isolate_beta}(b).

Repeating the above process, we can make sure that every $\alpha$ circle intersects some $\beta$ circle and every $\beta$ circle intersects some $\alpha$ circle.

\begin{figure}[htbp]
 \center{\includegraphics{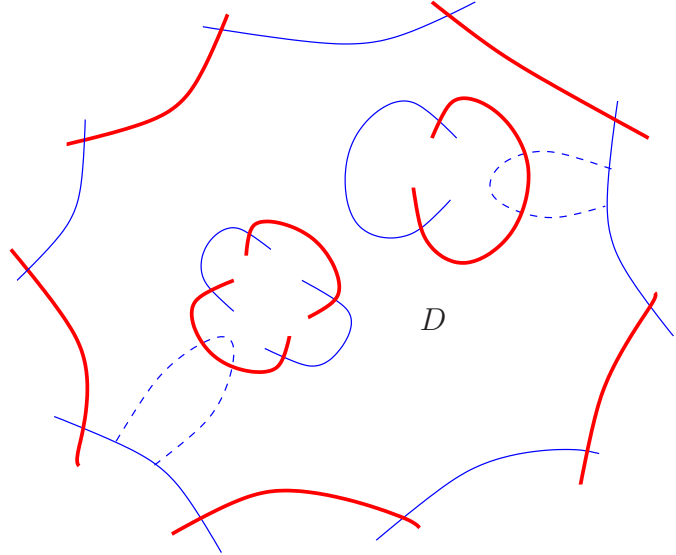}}
 \caption{\label{fig:kill_non_disk}
  {\bf Killing non-disk regions.}
  The dotted arcs indicate our finger moves. After our finger move, the region $D$ becomes a disk region. }
\end{figure}

Note that the complement of the $\alpha$ curves is a punctured sphere. Thus every region is a planar surface. A non-disk region $D$ has more than one boundary component. Every boundary component must contain both $\alpha$ and $\beta$ arcs since every $\alpha$ (resp. $\beta$) circle intersects some $\beta$ (resp. $\alpha$) circle. Then we make a finger move on the $\beta$ curve to reduce the number of boundary components of $D$ without generating other non-disk regions. See Figure \ref{fig:kill_non_disk} for this finger move operation. Repeating this process as many times as necessary, we will kill all the non-disk regions.

\subsection*{Step 2. Making all but one region bigons or squares.}\

We consider Heegaard diagrams with only disk regions. Note that our algorithm will not generate non-disk regions.

Let $D_0$ be the disk region containing the basepoint $w$. For any region $D$, pick an interior point $w^\prime\in D$ and define the {\it distance} of $D$, denoted by $d(D)$, to be the smallest number of intersection points between the $\beta$ curves and an arc connecting $w$ and $w^\prime$ in the complement of the $\alpha$ circles. For a $2n$-gon disk region $D$, define the {\it badness} of $D$ as $b(D)=max\{n-2,0\}$.

For a pointed Heegaard diagram $\calh$ with only disk regions, define the {\it distance} $d(\calh)$ of $\calh$ to be the largest distance of bad regions. Define the {\it distance $d$ complexity} of $\calh$ to be tuple
$$c_d(\calh)=\left(\sum_{i=1}^m b(D_i), -b(D_1), -b(D_2),\cdots, -b(D_m)\right),$$
where $D_1,\cdots,D_m$ are all the distance $d$ bad regions, ordered so that $b(D_1)\geq b(D_2)\geq \cdots \geq b(D_m)$. We call the first term the {\it total badness of distance $d$} of $\calh$, and denote it by $b_d(\calh)$. If there are no distance $d$ bad regions, then $c_d(\calh)=(0)$. We order the set of distance $d$ complexities lexicographically.

\begin{lemma} For a distance $d$ pointed Heegaard diagram $\calh$ with only disk regions, if $c_d(\calh)\neq(0)$, we can modify $\calh$ by isotopies and handleslides to get a new Heegaard diagram $\calh^\prime$ with only disk regions, satisfying $d(\calh^\prime)\leq d(\calh)$ and $c_d(\calh^\prime)<c_d(\calh)$.
\end{lemma}

\begin{proof} We order the bad regions of distance $d$ as in the definition of the distance $d$ complexity. Now we look at $D_m$. It is a $(2n)$-gon with $n\geq3$. Pick an adjacent region $D_*$ with distance $d-1$ having a common $\beta$ edge with $D_m$. Let $b_*$ be (one of) their common $\beta$ edge(s). We order the $\alpha$ edges of $D_m$ counterclockwise, and denote them by
$a_1,a_2,\cdots,a_n$ starting at $b_*$.

We try to make a finger move on $b_*$ into the $D_m$ and out of $D_m$ through $a_2$, as indicated in Figure \ref{fig:finger_move_start} when $D_m$ is an octagon. Our finger will separate $D_m$ into two parts, $D_{m,1}$ and $D_{m,2}$.
\begin{figure}[htbp]
 \center{\includegraphics{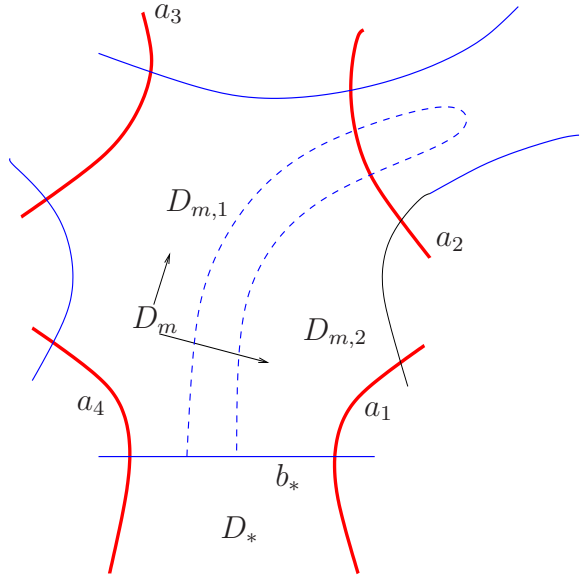}}
 \caption{\label{fig:finger_move_start}
  {\bf Starting our finger move.}
 }
\end{figure}

If we reach a square region of distance $\geq d$, we push up our finger outside the region via the opposite edge, as in Figure \ref{fig:finger_move_square}. Note that doing a finger move through regions of distance $\geq d$ does not change the distance of any of the bad regions, since they all have distance $\leq d$.
\begin{figure}[htbp]
 \center{
  \includegraphics{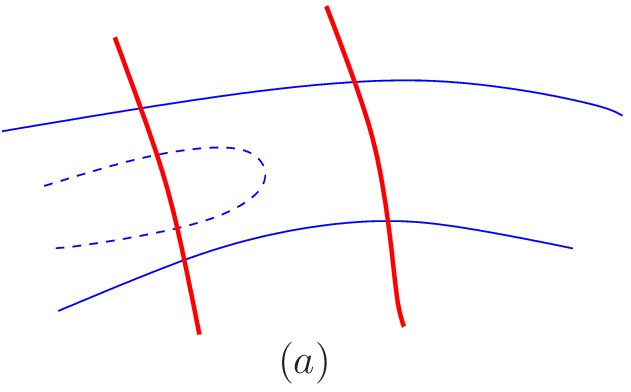}
  \hspace{30pt}
  \includegraphics{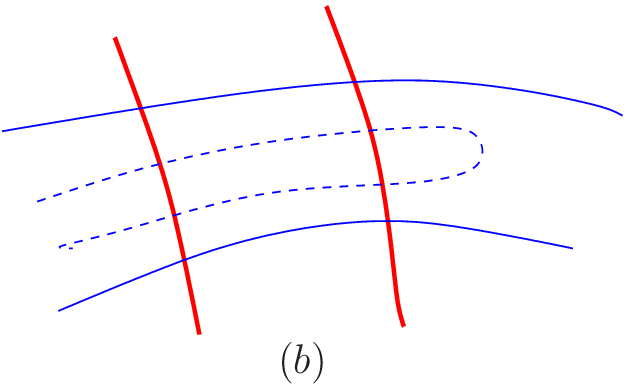}
 }
 \caption{\label{fig:finger_move_square}
  {\bf Moving across a square region.}
 }
\end{figure}

We continue to push up our finger as far as possible, until we reach one of the following:
\begin{enumerate}
 \item a bigon region.
 \item a region with distance $\leq d-1$.
 \item a bad region with distance $d$ other than $D_m$, i.e., $D_i$ with $i<m$.
 \item $D_m$.
\end{enumerate}
We will prove our lemma case by case.

\subsection*{Case 1. A bigon is reached.}\

Before we reach the bigon region, all regions in between are square regions with distance $\geq d$. After our finger moves inside a bigon region, our finger separates the bigon into a square and a new bigon, as in Figure \ref{fig:finger_move_bigon}.
\begin{figure}[htbp]
 \center{\includegraphics[width=2in,height=1.8in]{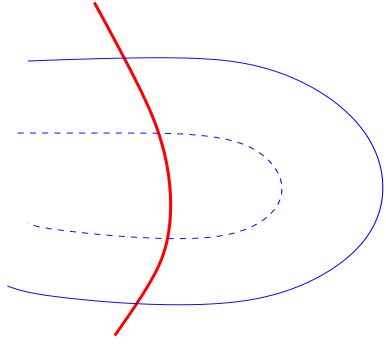}}
 \caption{\label{fig:finger_move_bigon}
  {\bf Case 1. A bigon is reached.}
 }
\end{figure}

Denote the new Heegaard diagram by $\calh^\prime$. We have $b(D_{m,1})=b(D_m)-1$. Since $D_{m,2}$ is a square and is good, we get $b_d(\calh^\prime)=b_d(\calh)-1$. Note that we will not increase the distance of any bad region since we do not pass through any region of distance $\leq d-1$ and all bad regions has distance $\leq d$. Hence $d(\calh^\prime)\leq d(\calh)$ and $c_d(\calh^\prime)<c_d(\calh)$.

\subsection*{Case 2. A smaller distance region is reached.}\

Let $D^\prime$ be the region with distance $<d$ we reached by our finger. Suppose $d(D^\prime)=d^\prime$. Let $\calh^\prime$ be the new Heegaard diagram. See Figure \ref{fig:finger_move_smaller_dist_region}. Note that $D^\prime$ might be a bigon, which could be covered in both Case 1 and Case 2.
\begin{figure}[htbp]
 \center{\includegraphics[width=3.5in,height=2.2in]{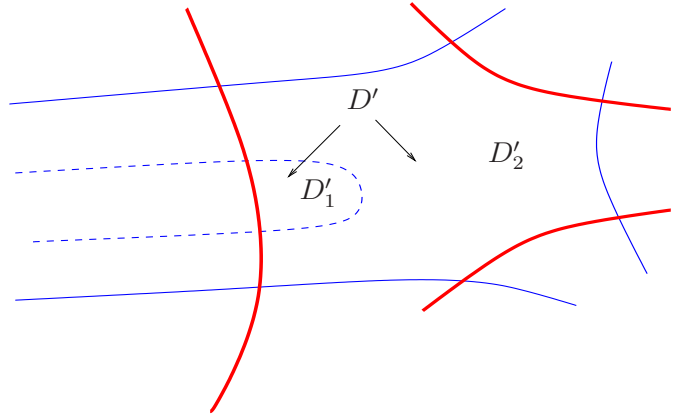}}
 \caption{\label{fig:finger_move_smaller_dist_region}
  {\bf Case 2. A smaller distance region is reached.}
 }
\end{figure}

We have $b(D_{m,1})=b(D_m)-1$ and $D_{m,2}$ is good. Our finger separates $D^\prime$ into a bigon region $D_1^\prime$ and the other part $D_2^\prime$. When $D^\prime$ is a square or a bad region, $D_2^\prime$ will be a bad region of distance $d^\prime<d$. We might have increased the distance $d^\prime$ complexity, but we have $d(\calh^\prime)\leq d(\calh)$ and $c_d(\calh^\prime)<c_d(\calh)$.

\subsection*{Case 3. Another distance $d$ bad region is reached.}\

In this case, we reach some distance $d$ bad region $D_i$ with $i<m$. See Figure \ref{fig:finger_move_brother_region} for an indication.
\begin{figure}[htbp]
 \center{\includegraphics[width=3.5in,height=2.2in]{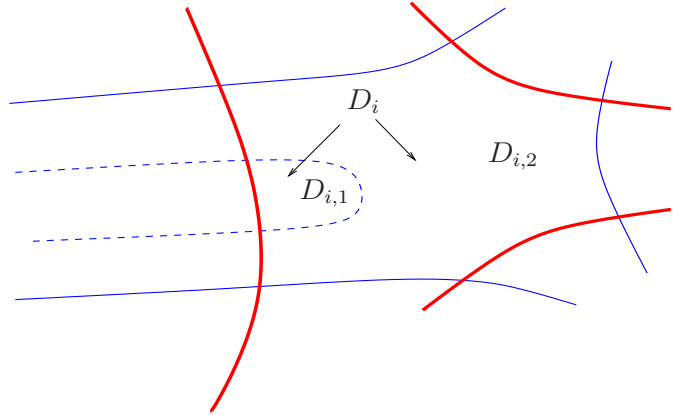}}
 \caption{\label{fig:finger_move_brother_region}
  {\bf Case 3. Another distance $d$ bad region is reached.}
 }
\end{figure}
Denote by $D_{i,1}$ and $D_{i,2}$ the two parts of $D_i$ separated by our finger. Then $D_{i,1}$ is good while $D_{i,2}$ is a bad region of distance $d$. We have $b(D_{i,2})=b(D_i)+1$ and $b(D_{m,1})=b(D_{m})-1$. Thus the total badness of distance $d$ remains the same. But we are decreasing the distance $d$ complexity since we are moving the badness from a later bad region to an earlier bad region. Hence for the new Heegaard diagram $\calh^\prime$, we have $d(\calh^\prime)=d(\calh)$ and $c_d(\calh^\prime)<c_d(\calh)$.

\subsection*{Case 4. Coming back to $D_m$.} This is the worst case and we need to pay more attention. We divide this case into two subcases, according to which edge the finger is coming back through.

\subsection*{Subcase 4.1. Coming back via an adjacent edge.}\

This subcase is indicated in Figure \ref{fig:finger_move_back_adjacent}.
\begin{figure}[htbp]
 \center{\includegraphics{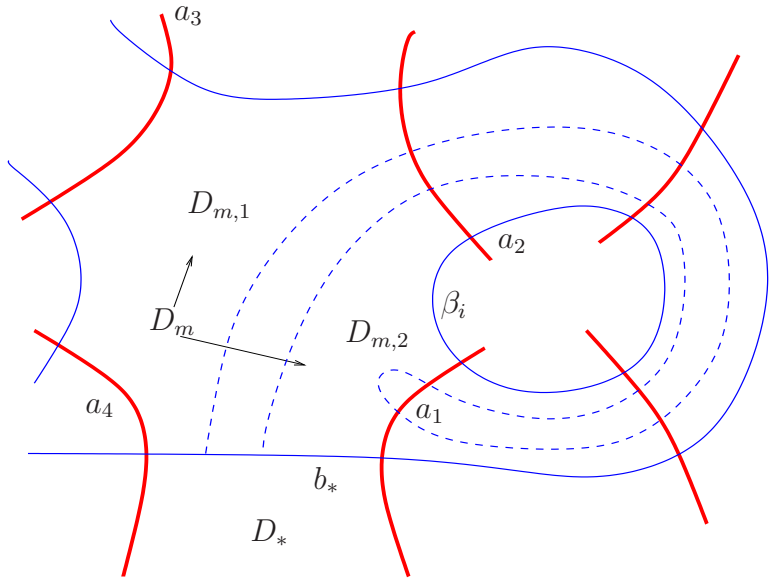}}
 \caption{\label{fig:finger_move_back_adjacent}
  {\bf Case 4.1 Coming back via an adjacent edge - finger move.} The finger is denoted by the dotted arc.
 }
\end{figure}
Without loss of generality, we assume the finger comes back via $a_1$. In this case, we see the full copy of some $\beta$ curve, say $\beta_i$, one the right side along our long finger. Suppose $b_*\subset \beta_j$. Note that $i\neq j$ since otherwise $b_*\subset \beta_i$ and we will reach either $D_m$ or $D_*$ at an earlier time. Now instead of doing the finger move, we handleslide $\beta_j$ over $\beta _i$. This is indicated in Figure \ref{fig:finger_move_back_adjacent_handleslide}.
\begin{figure}[htbp]
 \vspace{-10pt}
 \center{\includegraphics[width=4in,height=2.9in]{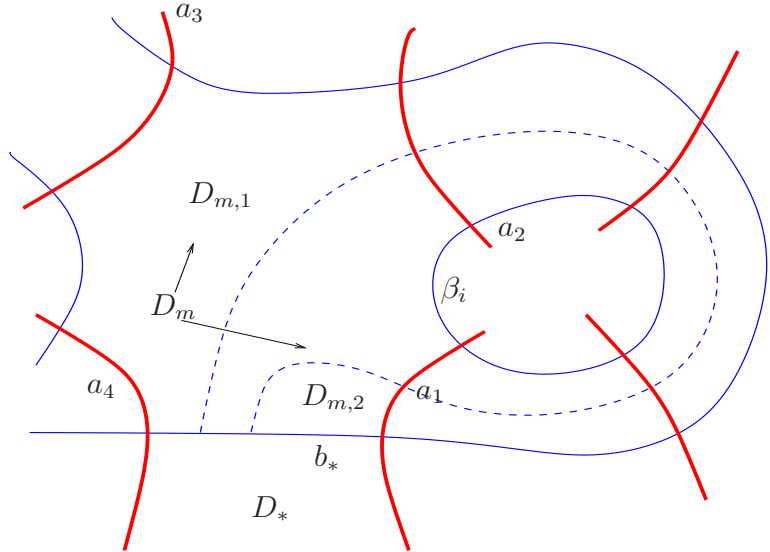}}
 \vspace{-10pt}
 \caption{\label{fig:finger_move_back_adjacent_handleslide}
  {\bf Case 4.1 \& 4.2 Coming back via an adjacent edge - handleslide.} The dotted arc denotes the $\beta$ curve after the handleslide.
 }
\end{figure}

Note that after the handle slides, we are not increasing the distance of any bad region. We have increased the badness of $D_*$, but it is a distance $d-1$ region. $D_{m,2}$ is a bigon region and $b(D_{m,1})=b(D_m)-1$. Thus for the new Heegaard diagram $\calh^\prime$ after the handleslide, the total badness of distance $d$ is decreased by 1. We have $d(\calh^\prime)\leq d(\calh)$ and $c_d(\calh^\prime)<c_d(\calh)$.

\subsection*{Subcase 4.2. Coming back via a non-adjacent edge.}\

If we return through $a_k$ with $3<k\leq n$, then, instead of the finger move through $a_2$, we do a finger move through $a_3$ (starting from $b_*$). If we reach one of the first three cases, we are decreasing the distance $d$ complexity by similar arguments as before.

Suppose instead that we come back to $D_m$, say via $a_i$. We claim that $3<i<k$.  Certainly we can not come back via $a_3$. The finger can not come back via $a_k$ since the chain of squares from $a_k$ is connected to $a_2$. If $i>k$ or $i<3$, we could close the cores the two fingers to get two simple closed curves $c_1$ and $c_2$, as indicated in Figure \ref{fig:no_crossing_finger}. 
\begin{figure}[htbp]
 \vspace{-18pt}
 \center{\includegraphics[width=280pt,height=270pt]{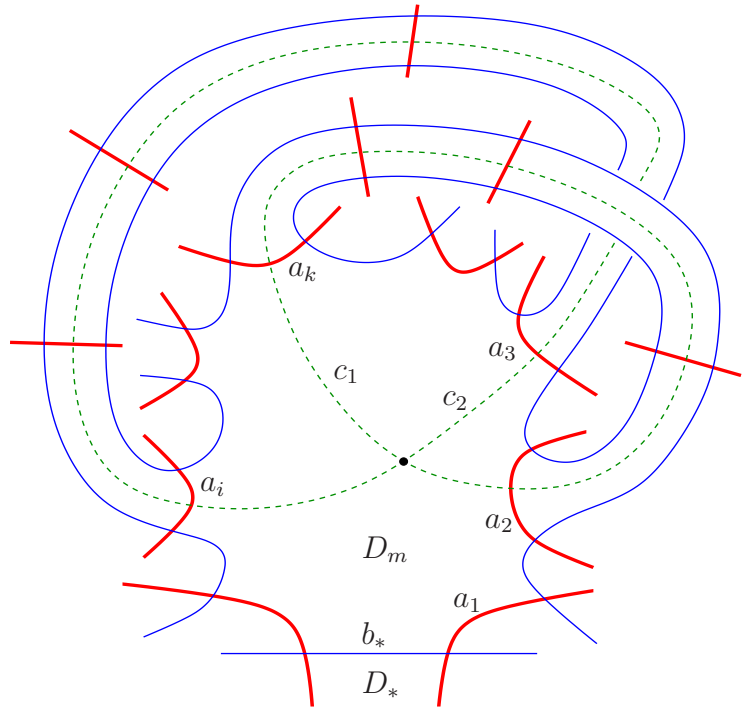}}
 \vspace{-10pt}
 \caption{\label{fig:no_crossing_finger}
  {\bf Case 4.2 There are no crossing fingers.} The fingers are not showed here. Instead, the two dotted arcs denote the cores of the two fingers.
 }
\end{figure}
Then $c_1$ and $c_2$ intersect transversely at exactly one point and they are in the complement of the $\beta$ curves. The complement of the $\beta$ curves is a punctured sphere. Attach disks to get a sphere. Then as homology classes, we get $[c_1]\cdot[c_2]=1$. But $H_1(S^2)\cong0$. This is a contradiction. Thus we must have $3<i<k$. (The argument of this claim was suggested by Dylan Thurston.)

Now, instead of the finger move through $a_3$, we do another finger move through $a_4$. Continuing the same arguments, we see that we either end up with a finger which does not come back, or we get some finger that starts at $a_j$ and comes back via $a_{j+1}$. If the finger does not come back, we reduce it to the previous cases and the lemma follows.

If there is a finger which starts at $a_j$ and comes back at $a_{j+1}$, we see a full $\beta$ circle. We do a handleslide similar to the one in Subcase 4.1. We have $b(D_{m,1})=\max\{n-j-1,0\}$ and $b(D_{m,2})=\max\{j-2,0\}$. We also have $b(D_{m,1})+b(D_{m,2})\leq n-3$. Thus for the new Heegaard diagram $\calh^\prime$ after the handleslide, the total badness of distance $d$ decreases. We have $d(\calh^\prime)\leq d(\calh)$ and $c_d(\calh^\prime)<c_d(\calh)$.

Thus we end the proof of our lemma.
\end{proof}

Repeat this process to make $c_d=(0)$. Repeating the whole process sufficiently many times will eventually kill all the bad regions other than $D_0$.

\subsection{Admissibility}

In this subsection, we show that our algorithm will not change the admissibility, that is, if we start with an admissible Heegaard diagram, then our algorithm ends with an admissible Heegaard diagram. There are two operations involved in our algorithm: isotopies and handleslides, and we will consider them one by one.

\begin{figure}[htbp]
 \center{
  \includegraphics{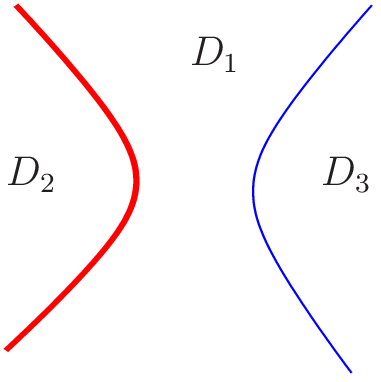}
  \hspace{30pt}
  \includegraphics{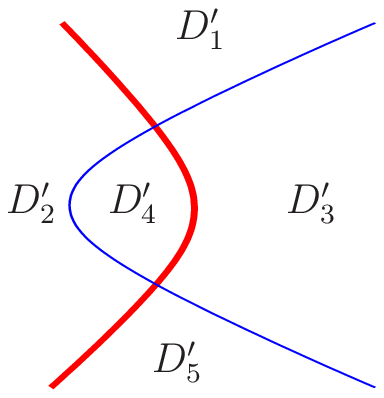}
 }
 \caption{\label{fig:admissible_isotopy}
  {\bf Isotopy of the $\beta$ curve.}
 }
\end{figure}

The isotopy is the operation in Figure \ref{fig:admissible_isotopy}. 
Let $\calh$ and $\calh^\prime$ be the Heegaard diagrams before and after the isotopy. Suppose $\calh$ is admissible. For a periodic domain in $\calh^\prime$
$$\phi^\prime=c_1D^\prime_1+c_2D^\prime_2+c_3D^\prime_3+c_4D^\prime_4+c_5D^\prime_5+\cdots$$
we have $c_2-c_1=c_4-c_3=c_2-c_5$ and $c_1-c_3=c_2-c_4=c_5-c_3$. Hence $c_1=c_5$ and $c_4=c_2+c_3-c_1$. Note that the regions are all the same except those in Figure \ref{fig:admissible_isotopy}. Therefore,
$$\phi=c_1 D_1+c_2 D_2 + c_3 D_3 + \cdots$$
is a periodic domain for $\calh$. Since $\calh$ is admissible, $\phi$ has both positive and negative coefficients, and so does $\phi^\prime$. Hence $\calh^\prime$ is admissible.

\begin{figure}[htbp]
 \center{
  \vspace{-30pt}
  \includegraphics{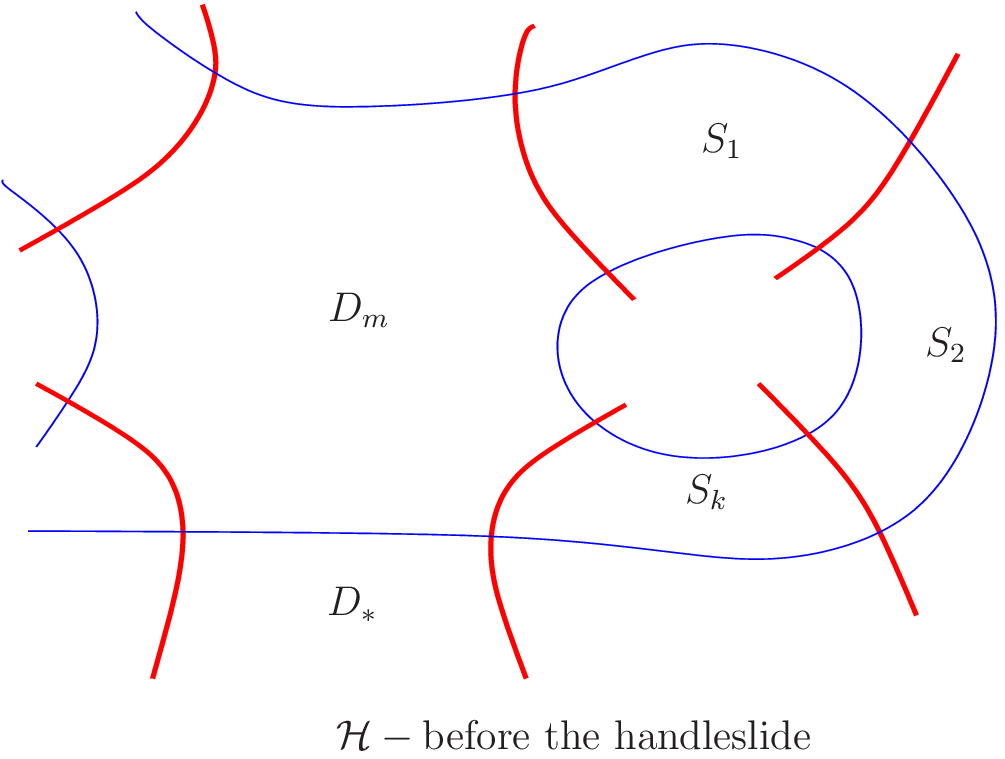}\\
  \vspace{-30pt}
  \includegraphics{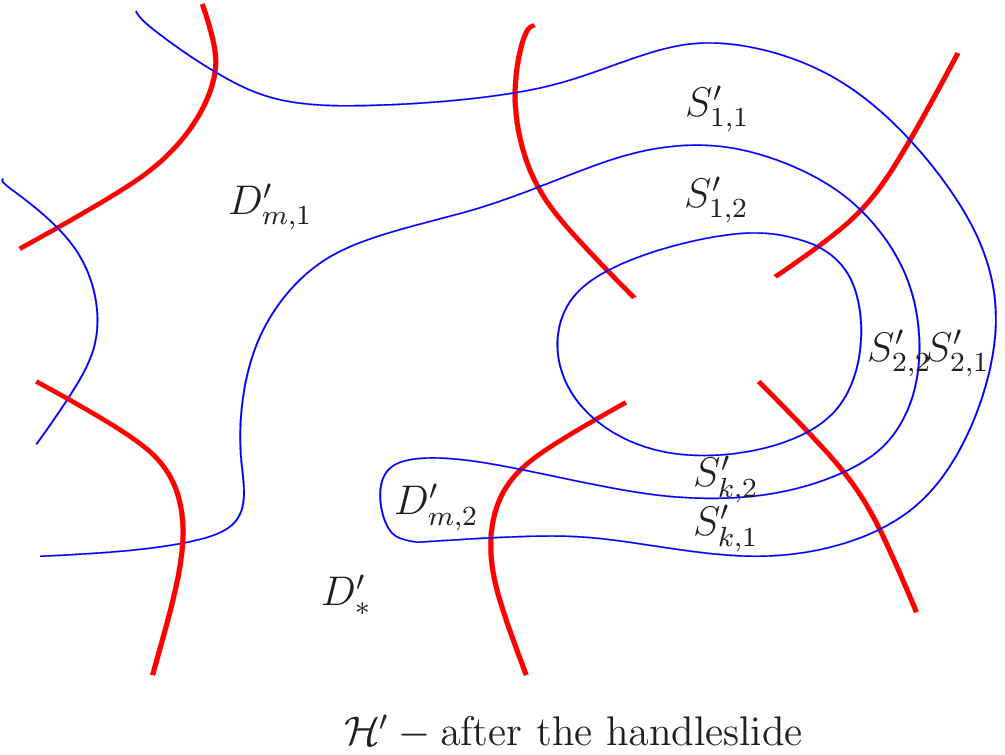}
  \vspace{-10pt}
 }
 \caption{\label{fig:admissible_handleslide}
  {\bf Handleslide of the $\beta$ curve.}
 }
\end{figure}

Our handleslide operation is indicated in Figure \ref{fig:admissible_handleslide}. 
Suppose $\calh$ is admissible. For a periodic domain in $\calh^\prime$
$$\phi^\prime=c_* D_*^\prime + c_1 D^\prime_{m,1} + c_{2} D^\prime_{m,2} + c_{1,1} S^\prime_{1,1} + c_{1,2} S^\prime_{1,2} +\cdots + c_{k,1} S^\prime_{k,1} + c_{k,2} S^\prime_{k,2} +\cdots$$
we get $c_{1}-c_*=c_{1,1}-c_{1,2}=\cdots=c_{k,1}-c_{k,2}=c_{2}-c_*$. Suppose $c_1-c_*=c_0$, then $c_{i,1}=c_{i,2}+c_0$ and $c_1=c_2$. Now
$$\phi=c_* D_* + c_1 D_m + c_{1,1} S_1+\cdots + c_{k,1} S_k+{\cdots}$$
is a periodic domain for $\calh$. Since $\calh$ is admissible, $\phi$ has both positive and negative coefficients. Hence $\phi^\prime$ has both positive and negative coefficients, so $\calh^\prime$ is admissible.

\vspace{10pt}

\noindent{\bf Remark.} In fact, it can be shown that nice Heegaard diagrams are always (weakly) admissible (\cite[Corollary 3.2]{LipManWang}).

We have similar conclusions for Heegaard diagrams with multiple basepoints. Our algorithm could be modified to get nice Heegaard diagrams in that case. Note that every region is connected to exactly one region containing some $w$ point in the complement of the $\alpha$ curves, so we can define the distance and hence the complexity in the same way, and thus our algorithm works as before. 

\begin{proof}[\bf Proof of Theorem \ref{thm:nicediagram}]
Starting with an admissible one-pointed Heegaard diagram, our algorithm described in Section 4.1 gives an admissible Heegaard diagram with only one bad region, the one containing the basepoint $w$. The algorithm can be modified for multiple basepoints as described above.
\end{proof}

\section{Examples}\label{section-examples}

In this section, we give two examples to demonstrate our algorithm. One is on knot Floer homology and the other is on the Heegaard Floer homology of three-manifolds.

\subsection{The Trefoil}

\begin{figure}[htbp]
 \center{\includegraphics{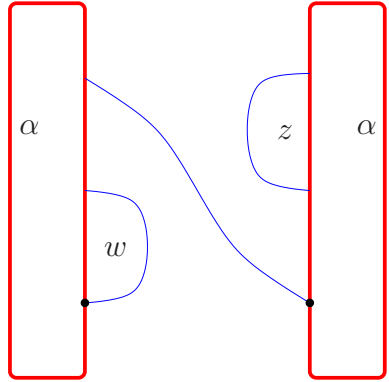}}
 \caption{\label{fig:trefoil_origin}
  {\bf A Heegaard diagram for the trefoil knot.} We make the convention that every two thick circles with the same $\alpha$ labels are identified so that the two dark points on them are identified.
 }
\end{figure}

We start with the Heegaard diagram of the trefoil knot in Figure \ref{fig:trefoil_origin}, where the two circles labeled by $\alpha$ are identified to get a genus one Heegaard diagram.

\begin{figure}[htbp]
 \center{\includegraphics[width=178pt]{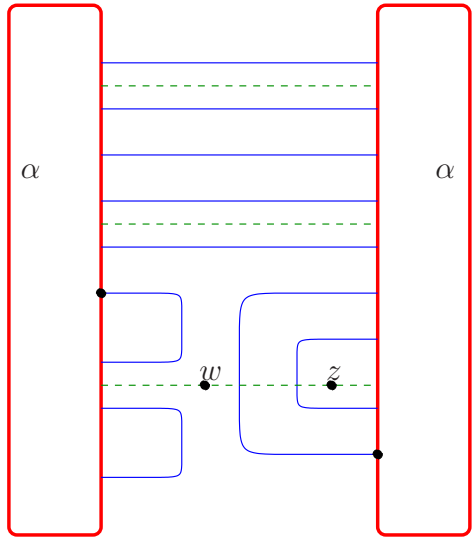}}
 \caption{\label{fig:trefoil_nice}
  {\bf A nice Heegaard diagram for the trefoil knot.} We use the same convention as in Figure \ref{fig:trefoil_origin}. The trefoil is given by the dotted curve.
 }
\end{figure}

After isotopy using the algorithm in Section \ref{section-algorithm}, we end up with the Heegaard diagram as in Figure \ref{fig:trefoil_nice}.
So we have nine generators. It is routine to find all boundary holomorphic disks and determine the Alexander and Maslov gradings of each generator.

\subsection{The Poincar{\'e} homology sphere $\Sigma(2,3,5)$} 
\begin{figure}[htbp]
 \center{\includegraphics{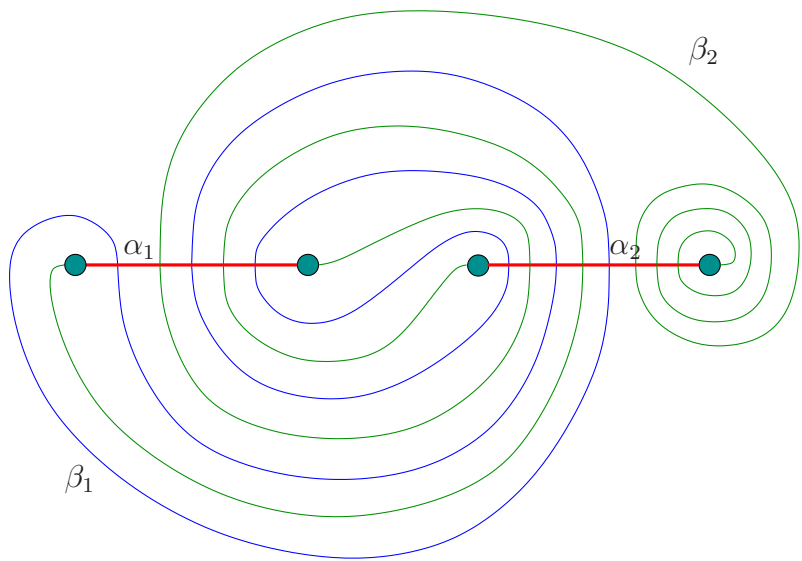}}
 \caption{\label{fig:poincare_from_trefoil}
  {\bf A Heegaard diagram for the Poincar{\'e} homology sphere.} The two darkly shaded circles on the left are the feet of one handle, and the two darkly shaded circles on the right are the feet of the other handle.
 }
\end{figure}

\begin{figure}[htbp]
 \center{\includegraphics[width=320pt]{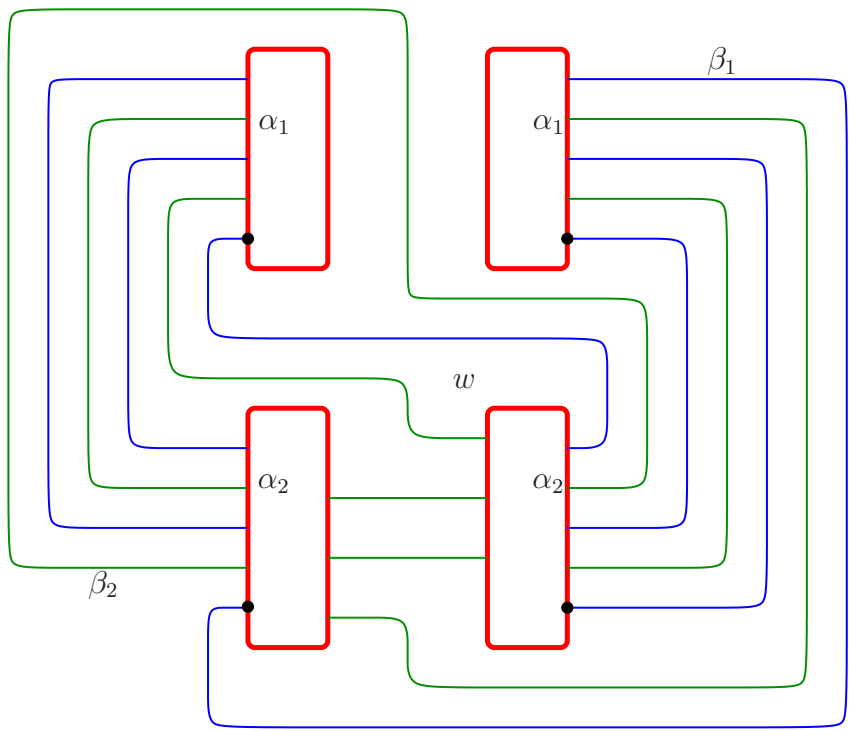}}
 \caption{\label{fig:poincare_origin}
  {\bf A Heegaard diagram for the Poincar{\'e} homology sphere.} We use the same convention as in Figure \ref{fig:trefoil_origin}.
 }
\end{figure}

We start with the Heegaard diagram of $\Sigma(2,3,5)$ in Figure \ref{fig:poincare_from_trefoil}, viewed as the +1 surgery on the right-handed trefoil knot. 
By cutting the Heegaard surface along the $\alpha$ circles, we get a planar presentation of the Heegaard diagram in Figure \ref{fig:poincare_origin}. 
It is easy to see that there are 21 generators for the chain complex. However, the authors do not know how to compute the differentials.

\begin{figure}[htbp]
 \center{\includegraphics{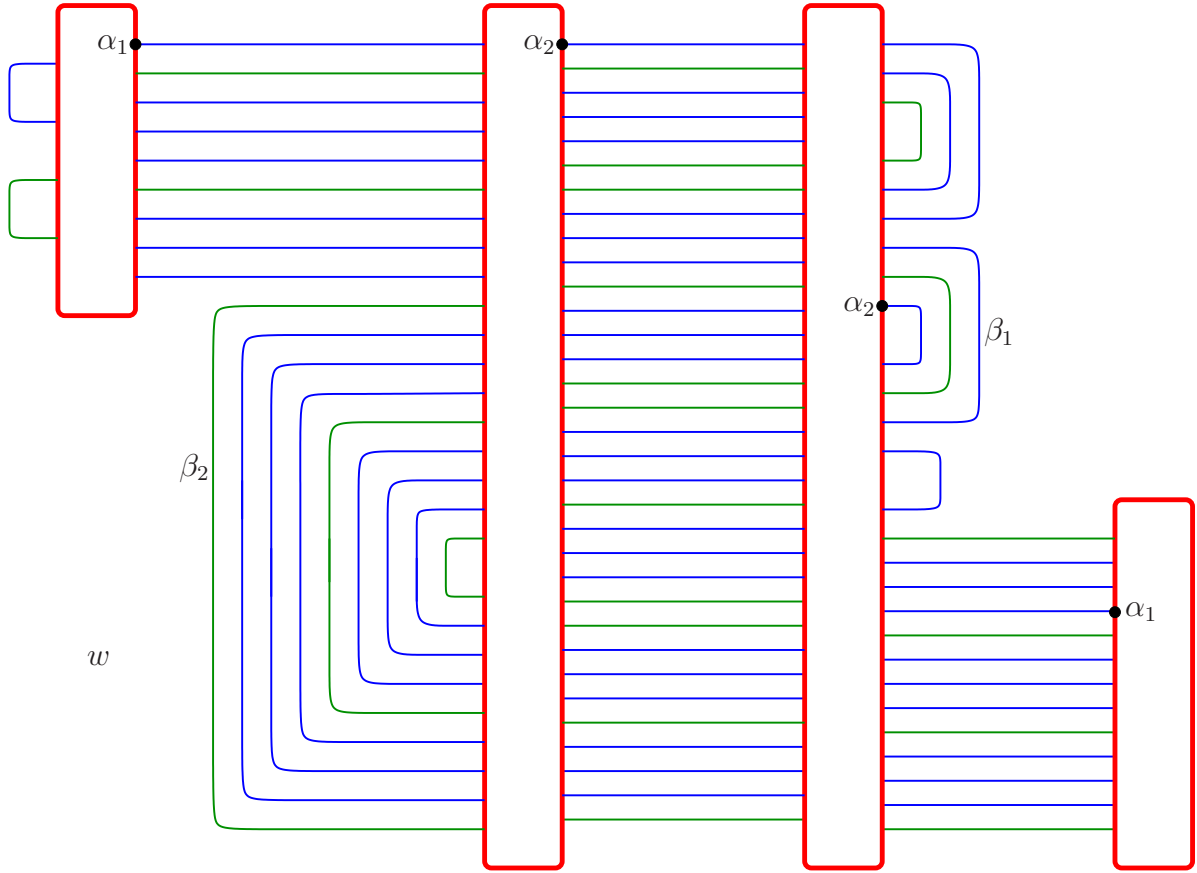}}
 \caption{\label{fig:poincare_nice}
  {\bf A nice Heegaard diagram for the Poincar{\'e} homology sphere.} We use the same convention as in Figure \ref{fig:trefoil_origin}. }
\end{figure}

After applying our algorithm, we get a nice Heegaard diagram as in Figure \ref{fig:poincare_nice}. 
There are 335 generators and 505 differentials for this diagram. We leave the actual computation using this diagram to the patient reader.

\end{document}